\documentclass[preprint,twocolumn,5p]{elsarticle}
  \def\appendixname{}
  \journal{Systems \& Control Letters}
\usepackage{mathtools,amsthm,amsmath,mathalfa,bbm,amssymb}
\usepackage{color}
\usepackage{comment}
\usepackage{tikz}
\usetikzlibrary{shapes,arrows}

\newcommand{\bbm}[1]{\left[\begin{matrix} #1 \end{matrix}\right]}
\newcommand{\sbm}[1]{\left[\begin{smallmatrix} #1
	\end{smallmatrix}\right]}

\newtheorem{thm}{Theorem}
\newtheorem{lem}{Lemma}
\newtheorem{prop}{Proposition}

\newtheorem{rem}{Remark}

\renewenvironment{proof}{\vspace{.1cm}\noindent{\sc Proof.}\hspace{0.10cm}\,\,}{$\hfill\Box$\vspace{.1cm}} 

\DeclareMathOperator*{\argmin}{arg\,min}

\newcommand{\dd}{\textrm{d}}

\let\norm\undefined % <-- "Undefine" \norm

\newcommand{\bracket}[1]{\left( #1 \right)}

\newcommand{\rline}{{\mathbb R}}

\newcommand{\rfb}[1]{\mbox{\rm
		(\eref{#1})}\ifx\undefined\stillediting\else:\fbox{$#1$}\fi}

\newcommand{\norm}[1]{\left\Vert#1\right\Vert}

\newcommand{\bluff}{{\hbox{\raise 15pt \hbox{\hskip 0.5pt}}}}
\newcommand{\etal}{\textit{et al }}
\newfont{\roma}{cmr10 scaled 1200}

\def\startmodif{\color{black}}

\begin{document}

%%============================================================================
\begin{frontmatter}
\title{Exponential Stabilization of Linear Systems using Nearest-Action Control with Countable Input Set}
\author[First,Second]{Muhammad Zaki Almuzakki\corref{cor1}}\ead{ m.z.almuzakki@universitaspertamina.ac.id}
\author[Second]{Bayu Jayawardhana}\ead{b.jayawardhana@rug.nl}
\author[Third]{Aneel Tanwani}\ead{aneel.tanwani@cnrs.fr}
\author[Second]{Antonis I. Vakis}\ead{a.vakis@rug.nl}

\cortext[cor1]{Corresponding author}
\address[First]{Department of Computer Science, Faculty of Science and Computer, Universitas Pertamina, Jakarta, Indonesia
}
\address[Second]{Engineering and Technology Institute Groningen, Faculty of Science and Engineering,  University of Groningen, The Netherlands
}
\address[Third]{Laboratory for Analyses and Architecture of Systems (LAAS) -- CNRS, Universit\'{e} de Toulouse, France
}

\tnotetext[fnGrant]{The work of M.Z.~Almuzakki is  supported by \emph{Lembaga Pengelola Dana Pendidikan Republik Indonesia} (LPDP-RI) under contract No. PRJ-851/LPDP.3/2016.}

\begin{abstract}
This paper studies stabilization of linear time-invariant (LTI) systems when control actions can only be realized in {\em finitely} many directions where it is possible to actuate uniformly or logarithmically extended positive scaling factors in each direction. Furthermore, a nearest-action selection approach is used to map the continuous measurements to a realizable action where we show that the approach satisfies a {\em weak} sector condition for multiple-input multiple-output (MIMO) systems. Using the notion of input-to-state
% practical 
stability, under some assumptions imposed on the transfer function of the system, we show that the closed-loop system converges to the target ball exponentially fast. Moreover, when logarithmic extension for the scaling factors is realizable, the closed-loop system is able to achieve asymptotic stability instead of only practical stability. Finally, we present an example of the application that confirms our analysis.
\end{abstract}


\end{frontmatter}
%%=========================================================================

\section{Introduction}

In some control systems applications, the actuator systems can only deliver constant input signals taken from a countable set of actuation points. 
For instance, the design of the Ocean Grazer's multi-piston-pump system provides only a finite number of constant actuation forces and limits the ability of the device to deliver arbitrary pumping forces; thus the control signals must be realizable based on the combination of multiple piston pumps 
\cite{Wei2018}. 

For solving such control problem where the control is taken from a countable input set, 
a nearest-action control (NAC) approach has recently been proposed and studied in \cite{JAYAWARDHANA2019460,Almuzakki2022}\footnote{In \cite{JAYAWARDHANA2019460,Almuzakki2022}, they use the term of nearest-neighbor control (NNC), which may be confused with the notion of ``neighbor'' in multi-agent systems. To avoid such confusion, the notion of ``action'' fits better than the former.} for a single agent and in \cite{Almuzakki2023} for multi-agent systems. It has been shown that for a given finite countable %, possibly finite, 
input set $\mathcal{U}$ satisfying some mild conditions, it is possible to render the closed-loop system practically stable by means of NAC. The ball, to which the state trajectories converge to, depends on the constellation of some elements in $\mathcal U$. Particularly, when the cardinality of $\mathcal U$ is minimal\footnote{By minimality, we mean the smallest number of elements of $\mathcal U$ that can be used to practically stabilized the systems using the nearest action control approach.}, the closed-loop system suffers in terms of control performance, e.g.\ linear convergence rate as confirmed by numerical results in \cite{Almuzakki2022,Almuzakki2023}.
In this article, we propose a generalization of our NAC approach which in particular yields exponential convergence.

In the literature, the feedback stabilization of dynamical systems using a finite or discrete input set is studied under the notion of quantized control, where the input space is partitioned into finite, or countably many, cells and a control action is assigned for each of these cells. Some classical references in this direction are \cite{Delchamps1990,BrocLibe00,tatikonda2000control,Elia2001,NairEvans04} and the reader can also consult a well-written survey \cite{NairFagn07}. In case of linear systems, it is observed that as long as the quantization error is bounded, the state trajectories are shown to converge exponentially to a ball around the origin whose size is determined by the maximum value of the error induced by the quantization. When we restrict our attention to initial conditions belonging to a compact set, it suffices to choose finitely many control actions to partition the corresponding compact set in the input space into finitely many cells with bounded quantization error. Techniques based on dynamic quantization have also been proposed to overcome the limitation of bounding the initial state so that the corresponding compact set (subject to quantization) becomes time-varying \cite{Libe03Aut}.

On the other hand, it is also possible to work with finitely many control actions with unbounded quantization errors which result in practical stability. In this regard, the authors in \cite{DEPERSIS2009602,dePersis2012,Jafarian2015} propose binary or ternary control systems. The underlying principal is to choose a control action which directs the state in the right direction without necessarily relying on keeping the quantization error small. These schemes lead to a {\em finite} countable input set defined on a regular grid and the exponential convergence property is no longer achieved, consistent with our previous findings based on NAC in \cite{JAYAWARDHANA2019460,Almuzakki2022,Almuzakki2023}. In fact, the NAC approach, proposed in \cite{JAYAWARDHANA2019460,Almuzakki2022}, describes a quantization based feedback using an optimization problem which induces a partition of the input space in finitely many cells. Similar questions, but without taking structure of dynamics into consideration, have been addressed in \cite{BullLibe06} which describes quantization regions using static optimization problems.

Motivated by the exponential convergence analysis in the aforementioned works on quantized control systems, we study in this paper the design and analysis of NAC systems with an infinite countable input set defined on an irregular grid. The approach builds on the basic ideas proposed in \cite{Almuzakki2022} where we focused on finding an input set of minimal cardinality $\mathcal{U}$ such that the dynamic system can be stabilized via a static output feedback. {\startmodif In this paper}, we make the transition from the finite input set to an infinite countable input set, which is designed based on an admissible minimal $\mathcal U$ that is extended to entire input space using either uniform or logarithmic partition. When the initial conditions are restricted to a compact set, it suffices to take $\mathcal U$ of finite cardinality and in this case the advantage compared to other approaches is that the number of symbols required for exponential convergence grows linearly with dimension of the input. This yields a considerable reduction in the size of input symbols required for exponential stabilization {\startmodif and we do so for systems with strictly positive real transfer function in the closed-loop. Similar to the approach proposed in this paper, Arvind \etal in \cite{Arvind2024} investigated the relaxation of the NAC method in  \cite{Almuzakki2022} by integrating switching mechanism akin to pulse width modulation (PWM) and by employing logarithmically-extended control points to achieve asymptotic stability of the closed-loop system.}

From the analysis viewpoint, the added flexibility in the input set allows us to get a bound on the quantization error and we can therefore, look at the closed-loop system as an asymptotically stable system with bounded quantization error. To analyze this closed-loop system, we use the approach based on {\em input-to-state stability}. Since the inception of this notion in the pioneering work \cite{Sontag89}, there has been significant interest to develop control techniques such that the resulting closed-loop is input-to-state stable with respect to unknown errors in the control action. This notion has been used in the quantized control \cite{Liberzon06}. For the quantization technique proposed in our work, we can think of the closed-loop system as the one driven by a sector bounded set-valued mapping plus the error due to quantization. Input-to-state stability for sector-bounded set-valued dynamics have been studied in \cite{Jayawardhana2009ISSincl,Jayawardhana2011} and here we use those results to develop the exponentially converging bounds on the norm of the state-trajectories in the closed-loop configuration.

The rest of this paper is organized as follows. In Section \ref{sec:prelims}, we provide some notations and preliminaries on the notion of practical stability and its relation to the ISS notion, on the nearest action control approach, and on the uniform and logarithmic quantizers. In Section \ref{sec:nnc_unifLoga}, we present our main results on the construction of infinite countable input set $\mathcal U$ and on the practical stability analysis of the resulting NAC systems. Furthermore, we give an illustrative example in Section \ref{sec:example} and 
conclusions in Section \ref{sec:conclusion}.

\section{Preliminaries and Problem Formulation}\label{sec:prelims}

{\bf Notation:} 
For a vector in $\rline^n$, or a matrix in $\rline^{m\times n}$, we denote the Euclidean norm and the corresponding induced norm by $\| \cdot \|$. A positive definite matrix $A\in \rline^{m\times m}$ is denoted by $A \succ 0$. 
For any point $c\in\mathbb R^n$, we define 
$\mathbb B^n_\epsilon(c) :=\{\xi\in \mathbb R^n | \|\xi-c\|\leq \epsilon\}$, and for simplicity, %we write 
$\mathbb B^n_\epsilon := \mathbb B^n_\epsilon(0)$. 
The inner product of 
$\mu,\nu\in \mathbb R^m$ is denoted by $\langle \mu, \nu \rangle$. For a given set $\mathcal{S}\subset\mathbb R^m$, and a vector $\mu \in \mathbb R^m$, $\langle \mu, \mathcal S \rangle := \{ \langle \mu, \nu \rangle \, \vert \, \nu \in \mathcal S \}$. For a countable set $\mathcal U\subset \rline^m$, its cardinality is denoted by $\text{card}(\mathcal U)$. The convex hull of vertices from a countable set $\mathcal U$ is denoted by $\text{conv}(\mathcal U)$. The interior of a set $S\subset\mathbb R^n$ is denoted by $\text{int}\bracket{S}$. A continuous function $\gamma:\mathbb R_{\geq 0}\to \mathbb R_{\geq 0}$ is of class $\mathcal K$ if it is continuous, strictly increasing, and $\gamma(0)=0$. We say that $\gamma:\mathbb R_{\geq 0}\to \mathbb R_{\geq 0}$ is of class $\mathcal K_\infty$ if $\gamma$ is of class $\mathcal K$ and $\lim_{s\to\infty}\gamma(s)=\infty$. A continuous function $\omega:\mathbb R_{\geq 0}\times\mathbb R_{\geq 0}\to \mathbb R_{\geq 0}$ is of class $\mathcal K\mathcal L$ if for each fixed $s$, $\omega(r,s)$ belongs to class $\mathcal K$, and for each fixed $r$, $\omega(r,s)$ is decreasing with respect to $s$ and is such that \ $\omega(r,s)\to 0$ as $s\to\infty$. The open right-half complex plane is denoted by $\mathbb C_+$. For a function $G:\mathbb C\to\mathbb C^{m\times m}$, we say that $G$ is strictly positive real if $G(s)+G^*(s)\succ 0,\ \forall s\in\mathbb C_+,\ s$ not a pole of $G$, where $G^*$ is the conjugate transpose of $G$. The space of $H^\infty$ is defined by $H^\infty:=\{G:\mathbb C\to\mathbb C^{m\times m}\ |\ G\text{ is holomorphic and }{\|G\|}_{H^\infty}:=\sup_{s\in\mathbb C_+}\|G(s)\|<\infty\}$ with $\|G(s)\|$ being the matrix norm induced by the $2$-norm on $\mathbb C^m$.

\subsection{Absolute stability, ISS \& practical stability}

{\startmodif
\begin{figure}
    \centering
    \tikzstyle{block} = [draw, rectangle, 
    minimum height=2em, minimum width=3em]
    \tikzstyle{sum} = [draw, circle, node distance=1cm]
    \tikzstyle{input} = [coordinate]
    \tikzstyle{output} = [coordinate]
    \tikzstyle{pinstyle} = [pin edge={to-,thin,black}]
    
    \begin{tikzpicture}[auto, node distance=2cm,>=latex']
        \node [input, name=input] {};
        \node [sum, right of=input] (sum) {};
        \node [block, right of=sum] (plant) {$(A,B,C)$};
        \node [block, below of=plant,
                node distance=1.cm] (controller) {$\Psi$};
    
        \draw [draw,->] (input) node[anchor=east] {$\Delta$} -- node[anchor=south west] {+} (sum);
        \draw [->] (sum) -- (plant);
        \draw [->] (plant) -- ([shift={(1cm,0cm)}]plant.east) |- (controller);
        \draw [->] (controller) -| (sum) node[anchor=north west, yshift=-0.2cm] {$-$};
    \end{tikzpicture}
    \caption{\startmodif 
    Block diagram of a Lur'e differential inclusion system, where a linear system with system matrices $(A,B,C)$ is feedback interconnected with a set-valued nonlinearity $\Psi$ and is subjected to a set-valued exogenous input signal $\Delta(t)$.}
    \label{fig:main_mimo_lure}
\end{figure}
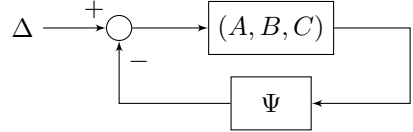
}

Consider the following Lur'e systems {\startmodif as depicted in Figure~\ref{fig:main_mimo_lure}}
\begin{equation}\label{sys:linear_nnc}
       \Sigma_{\rm lin} : \left\{ \begin{array}{rl} \dot x(t) & = Ax(t) + Bu(t) \\
        y(t) & = Cx(t) \\
        u(t) & \in \Delta(t)-\Psi(y(t)),
        \end{array} \right.
\end{equation}
where $x(t)\in\mathbb R^n$ is the state, $u(t),y(t)\in\mathbb R^m$ are the input and output, respectively, $\Delta : \mathbb R_{\ge 0} \to \mathbb R^m$ is a locally essentially bounded and a locally integrable function, 
the matrices $A,B,$ and $C$ are real matrices of suitable dimension, and $\Psi:\mathbb R^m \rightrightarrows \mathcal U$ with $\mathcal U\subseteq \mathbb R^m$ is a set-valued nonlinearity. The function $\Delta(t)$ embeds all exogeneous signals, possibly unknown, that can enter the nonsmooth system.  
When $\Delta \equiv 0$, the system \eqref{sys:linear_nnc} can be considered as a closed-loop interconnection of a linear differential inclusion system with system matrices $(A,B,C)$ and a set-valued nonlinearity $\Psi$. The transfer function of $\Sigma_{\rm lin}$ is given by $G(s)=C{(sI-A)}^{-1}B$.
{\startmodif The system \eqref{sys:linear_nnc} is said to be input-to-state stable (ISS) if there exist $\beta\in\mathcal K\mathcal L$ and $\rho\in\mathcal K_\infty$ such that, for all $x(0)\in\mathbb R^n$, the solution of \eqref{sys:linear_nnc} satisfies
\begin{equation}
    \|x(t)\|\le \beta(\|x(0)\|,t) + \rho({\|\Delta_{[0,t]}\|}_\infty),\ \forall t \ge 0.
\end{equation}
}

In \cite{Jayawardhana2009ISSincl,Jayawardhana2011}, the {\startmodif ISS} property of 
\eqref{sys:linear_nnc} has been established based on the system matrices $(A,B,C)$ and the sector bound condition of $\Psi$.
Although the ISS property of MIMO linear system \eqref{sys:linear_nnc} was established in \cite{Jayawardhana2009ISSincl,Sarkans2015}, the stability property requires a strong condition on the sector bound, namely $\langle k_1 y - v, k_2 y - v\rangle\le 0,\ \forall v\in \Psi(y),\ \forall y\in\mathbb R^m$ for some scalars $ k_1<k_2$. In this paper, we use the more familiar yet weaker version of the sector condition than before, that is $k_1{\|y\|}^2\le\langle v,y\rangle\le k_2{\|y\|}^2,\ \forall v\in \Psi(y),\ \forall y\in\mathbb R^m$ for some scalars $k_1<k_2$. One can immediately verify that the former condition implies the latter but not vice versa. Indeed, by multiplying the former sector condition with ${\|y\|}^2$ and since, by definition, we have ${\langle v,y\rangle}^2\le {\|v\|}^2 {\|y\|}^2$, it follows that
\begin{align*}
    0&\ge{\|y\|}^2\left\langle k_1 y - v, k_2 y - v\right\rangle + {\langle v,y\rangle}^2- {\|v\|}^2 {\|y\|}^2\\
    &= k_1 k_2 {\|y\|}^4 - (k_1+k_2)\langle v,y\rangle{\|y\|}^2{\langle v,y\rangle}^2\\
    &=\langle k_1 y - v, y\rangle \langle k_2 y - v, y\rangle.
\end{align*}
Since $k_1<k_2$, it must be that $\langle k_1 y - v, y\rangle\le0$ and $\langle k_2 y - v, y\rangle\ge 0$ which directly implies $k_1{\|y\|}^2\le\langle v,y\rangle\le k_2{\|y\|}^2$. For the converse, we can easily find an example where the latter sector condition does not imply the former. For example, by taking $k_1=1,\ k_2=2,\ y=\sbm{0\\2},$ and $v=\sbm{1\\2}$, it can be verified easily that the latter sector condition holds but the former does not.

\begin{rem}\label{rem:sector_cond}
    It is useful to note that the sector condition 
    \[
    k_1{\|y\|}^2\le\langle v,y\rangle\le k_2{\|y\|}^2,\ \forall v\in \Psi(y),\ y\in\mathbb R^m,\] 
    for some scalars $k_1<k_2$ implies
    \begin{equation}\label{ineq:sector_bound}
        \left\| v-\frac{k_1+k_2}{2}y\right\|\le\frac{k_2-k_1}{2c}\|y\|
    \end{equation}
    where $c\in(0,1]$ if $v\ne \frac{k_1+k_2}{2}y$ or $\| v-\frac{k_1+k_2}{2}y\|=0$, otherwise.
\end{rem}

The adaptation of the {\startmodif ISS} property of \eqref{sys:linear_nnc} in \cite[Theorem~3.4]{Jayawardhana2009ISSincl} is stated in the following theorem.

\begin{thm}\label{thm:iss_weak_sector}
    Consider the system $\Sigma_{\rm lin}$ in \eqref{sys:linear_nnc}. Suppose that the pair $(A,B)$ is controllable and $(A,C)$ is observable. 
    For the mapping $\Psi$, assume that there exist scalars $k_1<k_2$ such that for all 
    $v\in \Psi(y)$ 
    and for all $y\in\mathbb R^m$, it holds that $k_1{\|y\|}^2\le\langle v,y\rangle\le k_2{\|y\|}^2$. {In addition, assume that $G{(I+k_1 G)}^{-1}\in H^\infty$ and that $(I+k_2 G){(I+k_1 G)}^{-1}$ is strictly positive real.} Then every maximal solution {\startmodif is} forward complete and there exist positive constants $c_1, c_2,$ and $\varepsilon$ such that, for all $x(0)\in \mathbb R^n$,  
    every solution $x$ satisfies
    \begin{equation}\label{ineq:iss_thm1}
        \|x(t)\|\le c_1 e^{-\varepsilon t}\|x(0)\| + c_2 {\|\Delta_{[0,t]}\|}_\infty,\ \forall t\in\mathbb R_{>0}.
    \end{equation}
\end{thm}

The proof follows directly from the proof of \cite[Theorem~3.4]{Jayawardhana2009ISSincl}. The only difference is the part where we use Remark~\ref{rem:sector_cond} instead of \cite[Remark~3.1]{Jayawardhana2009ISSincl} in obtaining the result.  

Throughout the rest of this paper, we assume the following condition. 
\begin{itemize}
    \item[(A0)] For the system $\Sigma_{\rm lin}$ in \eqref{sys:linear_nnc}, the pair $(A,B)$ is controllable and $(A,C)$ is observable.
    
    \item[(A1)] For the system $\Sigma_{\rm lin}$ in \eqref{sys:linear_nnc} and the set-valued map ${\Psi}$, there exist scalars $k_1<k_2$ such that 
    \begin{equation}\label{ineq:sector_cond_weak}
        k_1{\|y\|}^2\le\langle v,y\rangle\le k_2{\|y\|}^2,\ \forall v\in  \Psi(y),\ \forall y\in\mathbb R^m;
    \end{equation}
    $G{(I+k_1 G)}^{-1}\in H^\infty$ and that $(I+k_2 G){(I+k_1 G)}^{-1}$ is strictly positive real.
\end{itemize}

Let {\startmodif $\rho\in\mathcal K_\infty$ and} $\omega$ be the minimum upper bound such that {\startmodif $\rho({\|\Delta(t)\|}_\infty)\le\omega$.}
% $c_2 {\|\Delta(t)\|}_\infty\le\omega$.
For the rest of the paper, we say that the system \eqref{sys:linear_nnc} is globally practically stable with respect to ball $\mathbb B_\omega$ ($\omega$-GPS) if \eqref{sys:linear_nnc} is {\startmodif ISS} with bias $\omega$, i.e.\ {\startmodif there exist $\beta\in\mathcal K\mathcal L$ so that} the global solution of $x$ satisfies
{\startmodif\[\|x(t)\|\le \beta(\|x(0)\|,t) + \omega,\ \forall t\in\mathbb R_{>0}.\]}
% \[\|x(t)\|\le c_1 e^{-\varepsilon t}\|x(0)\| + \omega,\ \forall t\in\mathbb R_{>0}.\]
If, in addition, the solution of \eqref{sys:linear_nnc} decays exponentially towards $\mathbb B_\omega$, we say that \eqref{sys:linear_nnc} is globally exponentially practically stable with respect to ball $\mathbb B_\omega$ ($\omega$-GEPS). Correspondingly, the system \eqref{sys:linear_nnc} is globally asymptotically stable (GAS) if it is $0$-GPS and it is globally exponentially stable (GES) if it is $0$-GEPS.

\subsection{Nearest Action Control}\label{sec:orig_NAC}

Given a set of {\startmodif nonzero} realizable actions {\startmodif $\mathcal U:=\{u_1, u_2, \dots, u_p\}$} satisfying $0\in{\rm int}({\rm conv}(\mathcal U))$, we assume the following assumption. 
\begin{itemize}
    \item[(A2)] For a given set {\startmodif $\mathcal U:= \{u_1, u_2, \ldots, u_p\}$}, there exists an index set $\mathcal I\subset \{1,\ldots,p\}$ such that the set $\mathcal V:=\{u_i\}_{i\in \mathcal I}\subset \mathcal U$ defines the vertices of a convex polytope satisfying, $0\in\text{int}\left(\text{conv}\left(\mathcal V\right)\right)$.
\end{itemize}

\begin{lem}[{ \cite[Lemma 1]{Almuzakki2022} }]\label{lemma:1}
Consider a discrete set $\mathcal U \subset \mathbb R^m$ that satisfies {\rm (A2)}. Then, there exists $\delta > 0$ such that
\begin{equation}\label{eq:bndVoronoi}
V_{\startmodif \mathcal U\cup \{0\}} (0) \subseteq \mathbb B_\delta,
\end{equation}
where $V_{\mathcal U}(s)$ is the Voronoi cell of $\mathcal U$ at a point $s\in\mathbb R^m$ defined by
\[V_\mathcal{U}(s) := \left\{x\in\mathbb R^m\ |\ \|x-s\|\leq\|x-v\|,\ \forall v\in\mathcal{U}\setminus\{s\} \right\}.\]
In other words, 
the following implication holds for each $\eta \in \mathbb R^m$
\begin{equation}\label{eq:defCv}
    \| \eta\| > \delta \Rightarrow \ \exists \ u_i \in {\startmodif \mathcal U\cup \{0\}} \text{ s.t. } \|\eta - u_i \|< \| \eta\|.
\end{equation}
\end{lem}

We define the nearest action map $\phi_{\mathcal U}:\mathbb R^m \rightrightarrows {\mathcal U}$ that maps any point $\eta\in\mathbb R^m$ to the nearest point $u\in{\mathcal U}$ as
\begin{equation}\label{phi_eq}
    \phi_{\mathcal U}(\eta):=\argmin_{v\in \mathcal{U}}\left\{{\|v-\eta\|}\right\}.
\end{equation}

\begin{lem}[{\cite{Almuzakki2022}}]\label{lem:phi_bndIneq}
    Consider a finite set {\startmodif $\mathcal U:=\{u_{1},u_{2},\dots,u_{p}\}$} satisfying {\rm (A2)} and the nearest action mapping $\phi_{\mathcal U}$ as in \eqref{phi_eq}. For a fixed $y\in\mathbb R^m$, let $\phi_{ \mathcal U}(-y)=\{u_j\}_{j\in\mathcal J}$ for some index set $\mathcal J\subset\{1,\dots,p\}$. Then the inequality
    \begin{equation}\label{eq:bndUpLowInnProd}
    -\|u_j\| \cdot \|y\| \le \langle u_j,y \rangle \le -\frac{1}{2} \|u_j\|^2
    \end{equation}
    holds for all $j\in\mathcal J$.
\end{lem}

In \cite{Almuzakki2022}, it has been shown that for a general class of passive nonlinear systems with proper storage function, large-time initial-state norm-observability assumption, and a given countable set of control actions {\startmodif $\mathcal U\cup \{0\}$}, the system can be practically stabilized by using the feedback law $u=\phi_{\startmodif \mathcal U\cup \{0\}}(-y)$. Note that for linear systems, the large-time norm-observability notion in nonlinear systems is equivalent to the usual observability notion for linear systems \cite[Remark 4]{hespanha2005nonlinear}. For the linear MIMO system, the unity output-feedback practical stabilization result in \cite{Almuzakki2022} can be expressed into the following proposition.
\begin{prop}\label{prop:1_aut}
    Consider the system $\Sigma_{\rm lin}$ in \eqref{sys:linear_nnc} with $\Delta \equiv 0$ satisfying {\rm (A0)} with $G(s)$ being strictly positive real, and a given finite set $\mathcal U$ satisfying {\rm (A2)} so that \eqref{eq:bndVoronoi} holds for some $\delta>0$. Let $\phi_{\mathcal U}$ be as defined in \eqref{phi_eq}. Let $\gamma: \mathbb{R}_{\ge 0} \to \mathbb{R}_{\ge 0}$ be defined as,%
    \footnote{The function $\gamma$ as described in \eqref{def_eq:gamma} is a natural consequence of the standard observability assumption for LTI systems. Interested readers may consult on \cite[Remark 1]{Almuzakki2022} for reference.}
\begin{equation}\label{def_eq:gamma}
        \gamma(\delta):=\delta\ \norm{{\bracket{W_\tau(t)}}^{-1}} \int_{t}^{t+\tau} \!\! \norm{e^{A^\top (s-t)}C^\top} \, \dd s, 
    \end{equation}
    where
    $W_\tau(t) = \int_{t}^{t+\tau} e^{A^\top (s-t)}C^\top Ce^{A (s-t)} \, \dd s$. 
    % {\startmodif and $\tau>0$. }%, and suppose 
    {\startmodif Suppose} that for a given $\omega>0${\startmodif, for each $t>0$, and for any $\tau>0$,} %{\startmodif and for all $\tau>0$}, 
    it holds that $\gamma(\delta)\le\omega$. Then, the closed-loop system with $u=\phi_{\startmodif\mathcal U\cup \{0\}}(-y)$ is $\omega$-GEPS.
\end{prop}

To prove the above proposition using Theorem~\ref{thm:iss_weak_sector}, we need to find suitable constants $k_1$ and $k_2$ that satisfy the weak sector condition \eqref{ineq:sector_cond_weak} with $\Psi(y)$ be replaced by $-\phi_{\startmodif\mathcal U\cup \{0\}}(-y)$ for all $y\in\mathbb R^m$. By Lemmas~\ref{lemma:1} and \ref{lem:phi_bndIneq}, for all $y\in\mathbb R^m\setminus\mathbb B_\delta$ and for all $w\in-\phi_{\startmodif\mathcal U\cup \{0\}}(-y)$, we have $w\ne0$ and 
$\frac{1}{2} \|w\|^2 \le \langle w,y \rangle \le \|w\| \cdot \|y\|$. The last inequality can be written as $k_1 \|y\|^2 \le \langle w,y \rangle \le k_2 \|y\|^2$ with 
\[k_1=\inf\limits_{\substack{w\in-\phi_{\startmodif\mathcal U\cup \{0\}}(-y)\\ y\in\mathbb R^m\setminus\mathbb B_\delta}}\frac{{\|w\|}^2}{2{\|y\|}^2}\text{ and }k_2=\sup\limits_{\substack{w\in-\phi_{\startmodif\mathcal U\cup \{0\}}(-y)\\ y\in\mathbb R^m\setminus\mathbb B_\delta}}\frac{\|w\|}{\|y\|}.\]
Then, we can conclude that the weak sector condition \eqref{ineq:sector_cond_weak} is satisfied outside the ball $\mathbb B_\delta$ with {\startmodif $\Psi(y)$} replaced by $-\phi_{\startmodif\mathcal U\cup \{0\}}(-y)$ for all $y\in\mathbb R^m\setminus\mathbb B_\delta$.
We also note that since the transfer function $G$ is already strictly positive real, it follows that the assumption (A1) is satisfied outside the ball $\mathbb B_\delta$.
{\startmodif Indeed, given that $G$ is strictly positive real, % and $0<k_1<k_2$, 
one can check that $G{(I+k_1 G)}^{-1}\in H^\infty$ since the system $(I+k_1G)^{-1}$ is a stable system for all $k_1>0$. %with input $u$ and output $y$ is a stable system. 
Additionally, it can be shown that the system $(I+k_2 G){(I+k_1 G)}^{-1}$ with $0<k_1<k_2$, with the corresponding new input $\tilde u=u+k_1 y$ and new output $\tilde y=k_2 y+u$ is also strictly positive real by the application of loop transformation to passive systems\footnote{For details on passivity and loop transformation, we refer interested reader to \cite{khalil2002nonlinear}.}, where the sector bound $[k_1,k_2]$ lies in the interval $[0,\infty)$.}

In order to apply Theorem~\ref{thm:iss_weak_sector}, we write the nearest action map $\phi_{\startmodif\mathcal U\cup \{0\}}$ as a linear combination of a mapping $\Psi_{\startmodif\mathcal U\cup \{0\}}(y)$ that satisfies the weak sector condition everywhere and an output-dependent disturbance term $\Delta_{\startmodif\mathcal U\cup \{0\}}(y)$ given by
\begin{equation}\label{eq:psi_phi1}
    \Psi_{\startmodif\mathcal U\cup \{0\}}(y):=\begin{cases}
        -\phi_{\startmodif\mathcal U\cup \{0\}}(-y), & y\in\mathbb R^m\setminus V_{\startmodif\mathcal U\cup \{0\}}(0),\\
        \frac{k_2+k_1}{2}y, & \text{otherwise.}
    \end{cases}
\end{equation}
and $\Delta_{\startmodif\mathcal U\cup \{0\}}(y):=\phi_{\startmodif\mathcal U\cup \{0\}}(-y)+\Psi_{\startmodif\mathcal U\cup \{0\}}(y)$, respectively. It follows that the mapping $\Psi_{\startmodif\mathcal U\cup \{0\}}$ in \eqref{eq:psi_phi1} satisfies assumption (A1) everywhere. Futhermore, by definition, for all $d\in\Delta_{\startmodif\mathcal U\cup \{0\}}$, we have that $\|d\|\le\delta$. Observe that the closed-loop of the linear system in $\Sigma_{\rm lin}$ with $u=\phi_{\startmodif\mathcal U\cup \{0\}}(y)$ is equivalent to $\Sigma_{\rm lin}$ with $\Delta$ replaced by $\Delta_{\startmodif\mathcal U\cup \{0\}}$ and $\Psi$ replaced by $\Psi_{\startmodif\mathcal U\cup \{0\}}$.
{\startmodif In addition, using the ideas similar to \cite[Remark 1]{Almuzakki2022}, we can write the constant $c_2$ in \eqref{ineq:iss_thm1} explicitly. In particular, since $\gamma(\delta)\le\omega$ for each $t>0$ and for any $\tau>0$, we can fix 
\[c_2=\sup\limits_{\substack{\forall t>0,\\\forall \tau>0}}\norm{{\bracket{W_\tau(t)}}^{-1}} \int_{t}^{t+\tau} \!\! \norm{e^{A^\top (s-t)}C^\top} \, \dd s\]
so that $\gamma(\delta)\le c_2\delta\le\omega$.
}
Finally, by direct application of Theorem~\ref{thm:iss_weak_sector} and since $\gamma(\delta)\le\omega$
we have that $\Sigma_{\rm lin}$ with $u=\phi_{\startmodif\mathcal U\cup \{0\}}(-y)=\Delta_{\startmodif\mathcal U\cup \{0\}}(y)-\Psi_{\startmodif\mathcal U\cup \{0\}}(y)$ is {\startmodif ISS} with bias {\startmodif $\omega$ , i.e.\ the closed-loop system is $\omega$-GEPS}.

\subsection{Uniform and Logarithmic Quantizers}

As briefly discussed in the Introduction, there are two standard types of quantization in literature, namely, the uniform and logarithmic quantizers. 
The range set of uniform quantizer is a regular grid and can be described by the set $\mathcal U_u^\lambda:=\{\pm k\lambda\ \mid\ k\in\mathbb Z_{\ge 0}\}$ with $\lambda>0$ be the desired uniform step size. One of the standard approach in the uniform quantization is the symmetric uniform quantizer given by
\begin{equation}\label{eq:uniform_standard}
    Q_u^\lambda(\eta) = \left\lfloor \frac{\eta}{\lambda} + \frac{1}{2} \right\rfloor \lambda.
\end{equation}
Observe that the symmetric uniform quantizer obeys the nearest-action rule, i.e.\ $Q_u^\lambda(\eta)=\phi_{\mathcal U_u^\lambda}(\eta)$. 

Similarly, the range set of logarithmic quantizer is a regular grid and is given by the set $\mathcal U_l^\lambda:=\{0\}\cup\{\pm \lambda^k\ \mid\ k\in\mathbb Z\}$ with $\lambda>1$ be the desired geometric step size. One example of the logarithmic quantizers is the mapping
\begin{equation}\label{eq:logarithmic_standard}
    Q_l^\lambda(\eta)=\begin{cases}
        0, & \eta=0\\
        \text{sign}(\eta) \lambda^{\left\lfloor \frac{1}{2}+\log_{\lambda}|\eta| \right\rfloor}, & \eta\ne 0.
    \end{cases}
\end{equation}
Note that both the uniform quantizer \eqref{eq:uniform_standard} and the logarithmic quantizer \eqref{eq:logarithmic_standard} are scalar functions. Typically, in the vectorized setting, the above quantizers are defined element-wise. 

\subsection{Nearest Action Control with Uniform and Logarithmic Points Extension}\label{sec:NAC_uniLoga_ext}

In this paper, we are interested in the output feedback stabilization of the system $\Sigma_{\rm lin}$ in \eqref{sys:linear_nnc}. In \cite{Almuzakki2022}, the nearest action control approach is presented where the input $u$ can only take values from a finite countable set ${\startmodif\mathcal U\cup \{0\}}$. 
As remarked above, in the standard multi-valued quantizers, the quantization takes place in each dimension of the input space leading to a regular grid of infinite countable input set $\mathcal U$. 

Instead of considering the regular grid obtained through element-wise quantization as before, we consider in this paper the extension of minimal ${\startmodif\mathcal U\cup \{0\}}$ studied in \cite{Almuzakki2022} to an infinite countable input set defined on the whole input space $\rline^m$. It is constructed by enlarging ${\startmodif\mathcal U\cup \{0\}}$ using each element $u_i\in\mathcal U$ as the vector to which infinitely new elements are generated uniformly or logarithmically. Without loss of generality, {for the rest of this paper, let us consider instead a finite countable set
\begin{equation}\label{eq:defCalU}
\begin{aligned}
& \mathcal{U} := \{u_1,\ldots,u_p \in \mathbb{R}^m : \| u_i\| = 1\} \\ 
& \quad \text{and } \mathcal U\bigcup\{0\} \text{ satisfies (A2).}
\end{aligned}
\end{equation}
Note that here we explicitly remove the zero (0) element from the set $\mathcal U$ used in {\startmodif \cite{Almuzakki2022,Almuzakki2023}.}
More formally, the uniformly-extended infinite countable set $\mathcal U_u^{\text{ext}}$ is given by
\begin{equation}\label{set:NNUQ}
    \mathcal U_u^{\text{ext},\lambda}:= \{k\lambda u\ \mid\ k\in\mathbb Z_{\ge 0},\ u\in\mathcal U\},
\end{equation}
where $\lambda>0$ is the desired uniform step size. Similarly, the logarithmically-extended infinite countable set $\mathcal U_l^{\text{ext}}$ is given by 
\begin{equation}\label{set:NNLQ}
    \mathcal U_l^{\text{ext},\lambda}:= \{0\}\bigcup\{\lambda^k u\ \mid\ k\in\mathbb Z,\ u\in\mathcal U\}
\end{equation}
where $\lambda>1$ is the desired geometric step size.
{
\begin{rem}
    The use of unit vectors $u_i$ in this extended set is to simplify the presentation of our main results and they are related to the characterization of the convergence ball. In general, we can consider any vectors of any length in the minimal countable set $\mathcal U$ to represent directions. Furthermore, if there exists at least one vector $u\in\mathcal U$ with $\|u\|\ne 1$, these vectors can be obtained by choosing a suitable step size $\lambda$ in $\mathcal U_u^{\text{ext},\lambda}$ or $\mathcal U_l^{\text{ext},\lambda}$.
\end{rem}
}

{\startmodif
\begin{rem}
    For the uniformly-extended infinite countable set $\mathcal U_u^{\text{ext},\lambda}$, the choice of the parameter $\lambda$ influences the size of the smallest ball centered at the origin that contains the Voronoi cell of the zero $(0)$ control input.
\end{rem}
}

By the definition of the nearest action map $\phi_{\mathcal U}$ in \eqref{phi_eq}, it is easy to see that for the uniformly distributed points, for all $z\in\mathbb R^m$, the nearest action mapping can be decomposed into
\begin{equation}\label{eq:nnc_phi_uniform}
    \phi_{\mathcal U_u^{\text{ext},\lambda}}(z)=\phi_{\mathcal U}(z)Q_u^\lambda\left(\langle z, \phi_{\mathcal U}(z) \rangle\right),
\end{equation}
with $Q_u^\lambda$ be the symmetric uniform quantizer in \eqref{eq:uniform_standard}. While the standard uniform quantizer $Q_u^\lambda$ defined in \eqref{eq:uniform_standard} obeys the nearest-action rule, the standard logarithmic quantizer $Q_l^\lambda$ defined in \eqref{eq:logarithmic_standard} does not. Instead, in the case of \eqref{eq:logarithmic_standard}, only the quantized exponent obeys the nearest-action rule, i.e.\ $-\frac{1}{2}\le \log_\lambda|\eta| - \left\lfloor \frac{1}{2}+\log_{\lambda}|\eta| \right\rfloor \le\frac{1}{2}$. The logarithmic quantizer \eqref{eq:logarithmic_standard} is therefore not suitable for the decomposition of $\phi_{\mathcal U_l^{\text{ext},\lambda}}$.

In order for the logarithmic quantizer to satisfy the nearest action rule, for any scalar $\eta\in\mathbb R_{>0}$ that is mapped to $\lambda^k,\ k\in\mathbb Z,\ \lambda\in\mathbb R_{>1}$, it must be that
\begin{align*}
    \frac{\lambda^{k-1}+\lambda^k}{2}\le \eta \le \frac{\lambda^{k}+\lambda^{k+1}}{2}.
\end{align*}
By inspecting the upper and lower bound of above inequality, we have that
\begin{align*}
    \eta \le \frac{\lambda^{k}+\lambda^{k+1}}{2} 
    \Leftrightarrow 
    \lambda^k &\ge \frac{2\eta}{\lambda+1}
    \Leftrightarrow
    k
    \ge \log_\lambda\left(\frac{2\eta}{\lambda+1}\right),
\end{align*}
and
\begin{align*}
    \eta \ge \frac{\lambda^{k-1}+\lambda^{k}}{2} 
    \Leftrightarrow 
    \lambda^k 
    \le \frac{2\lambda\eta}{\lambda+1}
    \Leftrightarrow 
    k 
    \le \log_\lambda\left(\frac{2\lambda\eta}{\lambda+1}\right),
\end{align*}
respectively. Using the above relationships and considering also the negative part of the input variable, we define the \emph{symmetric} logarithmic quantizer $Q_{sl}^\lambda:\mathbb R \to \mathbb R$ as
\begin{equation}\label{eq:logarithmic_symmetric}
    Q_{sl}^\lambda(\eta)=\begin{cases}
        0, & \eta=0\\
        \text{sign}(\eta) \lambda^{\left\lfloor \log_\lambda\left(\frac{2\lambda |\eta|}{\lambda+1}\right) \right\rfloor}, & \eta\ne 0.
        \end{cases}
\end{equation}
By using $Q_{sl}^\lambda$, for all $z\in\mathbb R^m$, we can decompose $\phi_{\mathcal U_l^{\text{ext},\lambda}}$ into 
\begin{equation}\label{eq:nnc_phi_log}
    \phi_{\mathcal U_l^{\text{ext},\lambda}}(z)=\phi_{\mathcal U}(z)Q_{sl}^\lambda\left(\langle z, \phi_{\mathcal U}(z) \rangle\right).
\end{equation}

Illustrations of the nearest action region obtained using the nearest action map $\phi_{\mathcal U_u^{\text{ext},\lambda}}$ in \eqref{eq:nnc_phi_uniform} or $\phi_{\mathcal U_l^{\text{ext},\lambda}}$ in \eqref{eq:nnc_phi_log} is shown in Fig.~\ref{fig:illustration_NAC_ext_region}. For the logarithmically-extended infinite countable set $\mathcal{U}_l^{\text{ext},\lambda}$, it can be seen from Fig.~\ref{fig:illustration_NAC_ext_region}(b) that the separating lines perpendicular to the direction of each $u_i\in\mathcal U$ (in blue) are equidistant to two black dots. Hence the name {\em symmetric} for \eqref{eq:logarithmic_symmetric}.

\begin{figure}[ht!]
    \centering
    \includegraphics[width=.49\textwidth]{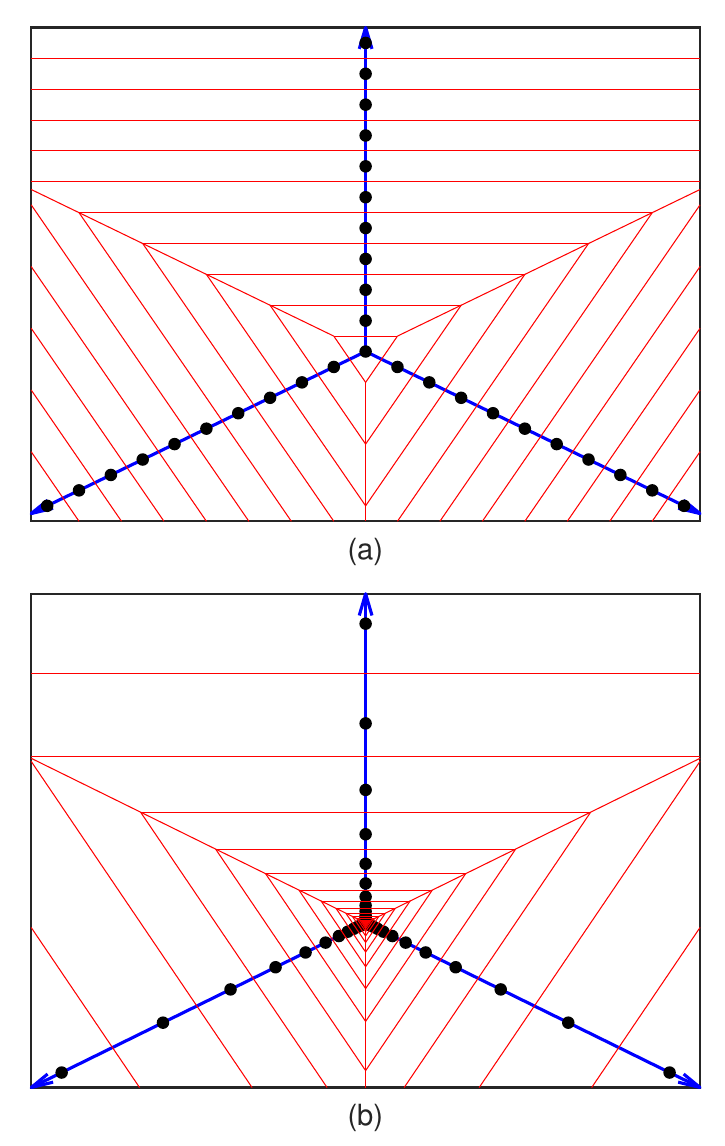}
    \caption{Illustration of the nearest action region of the (possibly infinite) realizable actions (represented by the black dots) set distributed (a) uniformly (as in \eqref{set:NNUQ}) or (b) logarithmically (as in \eqref{set:NNLQ}) along limited directions (blue arrows) in 2-dimensional input/output space. Here, the central coordinate is 0 (zero action). The region around each action (black dot) enclosed by the red lines represents the Voronoi cell of the respective action, i.e.\ all points in the enclosed area are mapped to their respective black dots by means of the nearest action map $\phi_{\mathcal U^{\text{ext}}}$.}
    \label{fig:illustration_NAC_ext_region}
\end{figure}

\section{Absolute Stability Analysis of The Nearest Action Control}\label{sec:nnc_unifLoga}

In this section, we present our main result on the output-feedback practical stabilization of the system $\Sigma_{\rm lin}$ described in \eqref{sys:linear_nnc} using the \emph{directional} nearest-action feedback control law $u=\phi_{\mathcal U}(-y)$ with $\phi_{\mathcal U}$ be as in \eqref{phi_eq}.

The set of realizable input considered in this section is the extended (possibly infinite) countable sets $\mathcal U^{\text{ext}}$ described in Section~\ref{sec:NAC_uniLoga_ext} along with their respective decomposable nearest action maps and scalar quantizers. In the general setting, we consider the extended set $\mathcal U^{\text{ext}}$ to be 
\begin{equation}\label{set:U_gen}
    \mathcal U^{\text{ext}}:=\{q u_i\ \mid\ u_i\in\mathcal U, q\in\mathcal Q\}
\end{equation}
where $\mathcal U$ is defined in \eqref{eq:defCalU} 
and $\{0,q_1, q_2,\dots\}=:\mathcal Q\subset\mathbb R_{\ge 0}$ with $0<q_1<q_2<\dots$ is a (possibly infinite) countable set containing non-negative scaling factors. %In addition, for 
Using the set $\mathcal Q$, we define a generic non-negative nearest action scalar quantizer $Q:\mathbb R_{\ge 0}\to\mathcal Q$ as follows
\begin{equation}\label{eq:general_scalar_quantizer}
    Q(\eta)=\phi_{\mathcal{Q}}(\eta):=\argmin\limits_{q\in\mathcal Q}\{\|q-\eta\|\}
\end{equation}
where $\eta\in\mathbb R_{\ge 0}$. 

Note that the quantizer $Q$ in \eqref{eq:general_scalar_quantizer} can be in any form of quantizer obeying the nearest action rule such as the symmetric uniform quantizer $Q_u^\lambda$ in \eqref{eq:uniform_standard}. Using \eqref{eq:general_scalar_quantizer}, for all $z\in\mathbb R^m$, the nearest action map $\phi_{\mathcal U^{\text{ext}}}$ as in \eqref{phi_eq} (with $\mathcal U$ be replaced by $\mathcal U^{\text{ext}})$ can be decomposed into
\begin{equation}\label{eq:general_NAC_ext}
    \phi_{\mathcal U^{\text{ext}}}(z)=\phi_{\mathcal U}(z)Q(\langle z,\phi_{\mathcal U}(z) \rangle).
\end{equation}
Then the following lemma for the setting in \eqref{set:U_gen}, \eqref{eq:general_scalar_quantizer}, and \eqref{eq:general_NAC_ext} is true.

\begin{lem}\label{lem:nnc_sectorBound}
    If for all $\eta\in[\frac{1}{2}q_{1},\infty)$, there exists $\kappa_1,\kappa_2\in\mathbb R_{>0}$ with $\kappa_1<\kappa_2$ so that the scalar quantizer $Q$ satisfies the sector bound
\begin{equation}\label{ineq:sector_bound_general}
        \kappa_1 \eta^2\le \eta Q(\eta) \le \kappa_2 \eta^2,
    \end{equation}
    then for all $z\in\mathbb R^n$ and $n\in\mathbb N$ satisfying $\langle z,\phi_{\mathcal{U}}(z) \rangle\ge \frac{1}{2}q_{1}$, the following inequality holds for some $\alpha\in(0,1]$,    \begin{equation}\label{ineq:nnc_sector_bound}
        \alpha\kappa_1 {\|z\|}^2\le \langle z,\phi_{\mathcal{U}^{\text{ext}}}(z)\rangle \le \kappa_2 {\|z\|}^2.
    \end{equation}
\end{lem}

\begin{proof}
    We first note that the nearest action map $\phi_{\mathcal U^{\text{ext}}}(z)$ is exactly the direction pointed by the mapping $\phi_{\mathcal U}(z)$ multiplied by the positive scalar obtained from quantizing the scalar projection of $z$ in the same direction using the scalar quantizer $Q$. We also note that since for all $u\in\mathcal U$, $\|u\|=1$, by the definition of $\phi_{\mathcal U}$ we have that $\frac{1}{2}\le\langle z, \phi_{\mathcal U}(z)\rangle \le {\|z\|}$. This means that there exists a lower bound $\alpha\in(0,1]$ so that the following holds for all $z\in\rline^n$
    \begin{equation}\label{ineq:nnc_innProd_zQU1}
        \alpha{\|z\|}\le\langle z, \phi_{\mathcal U}(z)\rangle \le {\|z\|}.
    \end{equation}
    By taking the upper bound in \eqref{ineq:sector_bound_general}, the decomposed nearest action map \eqref{eq:general_NAC_ext}, and using the inequality \eqref{ineq:nnc_innProd_zQU1}, it follows that for all $z\in\mathbb R^n$ satisfying $\langle z,\phi_{\mathcal{U}}(z) \rangle\ge \frac{1}{2}q_{1}$ we have
    \begin{align*}
        & \langle z, \phi_{\mathcal U^{\text{ext}}}(z) \rangle
        = \langle z, \phi_{\mathcal U}(z) Q\left(\langle z, \phi_{\mathcal U}(z) \rangle\right) \rangle\\
        &= \langle z, \phi_{\mathcal U}(z) \rangle Q\left(\langle z, \phi_{\mathcal U}(z) \rangle\right) 
        \le \kappa_2 {\langle z, \phi_{\mathcal U}(z) \rangle}^2 
        \le \kappa_2 {\|z\|}^2.
    \end{align*}
    Similarly, for the lower bound in \eqref{ineq:sector_bound_general} we have that
    \begin{align*}
        \langle z, &\phi_{\mathcal U^{\text{ext}}}(z) \rangle
        = \langle z, \phi_{\mathcal U}(z) \rangle Q\left(\langle z, \phi_{\mathcal U}(z) \rangle\right) \\
        &\ge \kappa_1 {\langle z, \phi_{\mathcal U}(z) \rangle}^2
        \ge \alpha\kappa_1 {\|z\|}^2.
    \end{align*}
    Therefore, the inequality \eqref{ineq:nnc_sector_bound} holds.
\end{proof}

Using the weak sector bound result on $\phi_{\mathcal U^{\text{ext}}}$ stated in Lemma~\ref{lem:nnc_sectorBound}, we can then analyze the stability property of the system $\Sigma_{\rm lin}$ described by \eqref{sys:linear_nnc} using the notion of {\startmodif ISS} as presented in the following proposition.

\begin{prop}\label{prop:main1}
    Consider the system $\Sigma_{\rm lin}$ in \eqref{sys:linear_nnc} satisfying {\rm (A0)} and a discrete set $\mathcal U^{\text{ext}}$ as in \eqref{set:U_gen} constructed from a finite countable set of unit vectors
    $\mathcal U$ that together with $\{0\}$ 
satisfies {\rm (A2)} so that \eqref{eq:bndVoronoi} holds for some $\delta>0$. Let $\phi_{\mathcal U^{\text{ext}}}$ be as in \eqref{eq:general_NAC_ext}. Suppose that $\phi_{\mathcal U^{\text{ext}}}$ satisfies \eqref{ineq:nnc_sector_bound} for some $\alpha\in(0,1]$ and $0<\kappa_1<\kappa_2$. In addition, assume that $G{(I+\alpha\kappa_1 G)}^{-1}\in H^\infty$ and that $(I+\kappa_2 G){(I+\alpha\kappa_1 G)}^{-1}$ is strictly positive real. Then the closed-loop system with $u=\phi_{\mathcal U^{\text{ext}}}(-y)$ is $\omega$-GEPS with $\omega=c_2 \delta q_1$ for some $c_2>0$.
\end{prop}

\begin{proof}
We first note that since by the definition of $\phi_{\mathcal U^{\text{ext}}}$ in \eqref{eq:general_NAC_ext} and according to Lemma~\ref{lemma:1}, for the set $\mathcal U^{\text{ext}}$ we have that {\startmodif $V_{\mathcal U^{\text{ext}}}(0)\subseteq{\mathbb B}_{\delta q_1}$}. In order to successfully apply Theorem 1, and since $\phi_{\mathcal U^{\text{ext}}}$ satisfies assumption (A1) only outside the ball ${\mathbb B}_{\delta q_1}$, we define $\Psi_{\mathcal U^{\text{ext}}}$ as in \eqref{eq:psi_phi1} with $\mathcal{U}$ replaced by $\mathcal U^{\text{ext}}$, $k_1=\alpha\kappa_1$, and $k_2=\kappa_2$, so that the sector condition is satisfied everywhere by means of $\Psi_{\mathcal U^{\text{ext}}}$. Following the proof of Proposition~\ref{prop:1_aut}, we can apply the result in Theorem~\ref{thm:iss_weak_sector} to conclude that the system $\Sigma_{\rm lin}$ is $\omega$-GEPS with $\omega=c_2 \delta q_1$ for some $c_2>0$.
\end{proof}

\subsection{Practical Stabilization with Uniformly-Extended Actions}

\begin{lem}\label{lem:uniform_sectorBound}
    Let $Q_u^\lambda:\mathbb R \to \mathbb R$ be the uniform quantizer defined in \eqref{eq:uniform_standard} along with a given desired stepsize $\lambda>0$. For all $\sigma \ge \frac{\lambda}{2}$, it holds that  %it holds that
\begin{equation}\label{ineq:uniform_quantizer_nnc}
        \left( 1-\frac{\lambda}{2\sigma} \right) \eta^2 \le \eta Q_u^\lambda(\eta) \le \left( 1+\frac{\lambda}{2\sigma} \right) \eta^2,
    \end{equation}    
    for all $\eta\ge\sigma$.
\end{lem}

\begin{proof}
    By the definition of $Q_u^\lambda(\eta)$ in \eqref{eq:uniform_standard}, the difference between $\eta$ and $Q_u^\lambda(\eta)$ satisfies
    \[-\frac{\lambda}{2}\le \eta-Q_u^\lambda(\eta)\le\frac{\lambda}{2}.\]
    Taking the upper bound of above inequality, it follows that
    \begin{align*}
        \eta-Q_u^\lambda(\eta)\le\frac{\lambda}{2} &&\iff&&
        Q_u^\lambda(\eta)&\ge \eta-\frac{\lambda}{2}=\left( 1-\frac{\lambda}{2\eta}\right)\eta\\
        &&&& &\ge\left( 1-\frac{\lambda}{2\sigma}\right)\eta.
    \end{align*}
    Similarly, using the lower bound $-\frac{\lambda}{2}\le \eta-Q_u^\lambda(\eta)$, it follows that
    \begin{align*}
        \eta-Q_u^\lambda(\eta)\ge-\frac{\lambda}{2} &&\iff&&
        Q_u^\lambda(\eta)&\le \eta+\frac{\lambda}{2}=\left( 1+\frac{\lambda}{2\eta}\right)\eta\\
        &&&& &\le\left( 1+\frac{\lambda}{2\sigma}\right)\eta.
    \end{align*}
    Finally, by combining the upper and lower bounds of $Q_u^\lambda(\eta)$ and multiplying all sides of the combined inequality by $\eta$ we have the inequality \eqref{ineq:uniform_quantizer_nnc}.
\end{proof}

\begin{prop}\label{cor:nnc_uniform}
    Consider the system $\Sigma_{\rm lin}$ described by \eqref{sys:linear_nnc} satisfying {\rm (A0)} and a discrete set $\mathcal U_u^{\text{ext},\lambda}$ as in \eqref{set:NNUQ} 
    constructed from a finite countable set of unit vectors $\mathcal U$ that together with $\{0\}$ 
satisfies {\rm (A2)} so that \eqref{eq:bndVoronoi} holds for some $\delta>0$. 
Let $\phi_{\mathcal U_u^{\text{ext},\lambda}}$ and $Q_u^\lambda$ be as in \eqref{eq:nnc_phi_uniform} and \eqref{eq:uniform_standard}, respectively, along with a given desired stepsize $\lambda>0$. 
Suppose that $\phi_{\mathcal U_u^{\text{ext},\lambda}}$ satisfies \eqref{ineq:nnc_sector_bound} for some $\alpha\in(0,1]$, $\kappa_1=\left(1-\frac{1}{2\delta}\right)$, and $\kappa_2=\left(1+\frac{1}{2\delta}\right)$. In addition, assume that $G{(I+\alpha\kappa_1 G)}^{-1}\in H^\infty$ and that $(I+\kappa_2 G){(I+\alpha \kappa_1 G)}^{-1}$ is strictly positive real. 
Then the closed-loop system with $u = \phi_{\mathcal U_u^{\text{ext},\lambda}}(-y)$ is $\omega$-GEPS with $\omega=c_2 \delta \lambda$ for some $c_2>0$.
\end{prop}

\begin{proof}
    We first observe that with regards to the nearest action selection approach where $\delta\ge\frac{1}{2}$ is the minimum upper bound satisfying $\mathcal U\subseteq\mathbb B_\delta$, the condition $\frac{\lambda}{2}\le\sigma$ in Lemma~\ref{lem:uniform_sectorBound} is satisfied with $\sigma=\lambda\delta$. Next, 
    for all $\eta\ge\lambda\delta$, we have that
    \[\left( 1-\frac{1}{2\delta} \right) \eta^2 \le \eta Q_u^\lambda(\eta) \le \left( 1+\frac{1}{2\delta} \right) \eta^2.\]
    Using the result in Lemma~\ref{lem:nnc_sectorBound}, it follows that for all {\startmodif $z\in\mathbb R^n$}, $n\in\mathbb N$, satisfying $\langle z,\phi_{\mathcal{U}}(z) \rangle\ge \frac{\delta}{2}$, there exists $\alpha\in(0,1]$ so that the mapping $\phi_{\mathcal{U}_u}(z)$ satisfies the sector condition
    \begin{equation}\label{eq:phi_unif}
        \alpha\left( 1-\frac{1}{2\delta} \right) {\|z\|}^2\le \langle z,\phi_{\mathcal{U}_u^{\text{ext},\lambda}}(z)\rangle \le \left( 1+\frac{1}{2\delta} \right) {\|z\|}^2.
    \end{equation}
    Finally, we use the result in Proposition~\ref{prop:main1} to complete the proof.
\end{proof}

\subsection{Global Exponential Stabilization with Logarithmically-Extended Actions}

We next present the case where the system can only realize actions in the direction contained in $\mathcal U$ with logarithmically distributed positive scaling factors. For this purpose, we have that the logarithmic quantizer $Q_{sl}^\lambda(\eta)$ defined in \eqref{eq:logarithmic_symmetric} satisfies the sector condition in (A1) as shown in the following lemma.

\begin{lem}\label{lem:symLogarithmic_sectorBound}
    Let $Q_{sl}^\lambda:\mathbb R \to \mathbb R$ be the symmetric logarithmic quantizer defined in \eqref{eq:logarithmic_symmetric} along with a given desired stepsize $\lambda>1$. Then \begin{equation}\label{ineq:symLogarithmic_quantizer_nnc}
        \left( \frac{2}{\lambda+1} \right) \eta^2 \le \eta Q_{sl}^\lambda(\eta) \le \left( \frac{2\lambda}{\lambda+1} \right) \eta^2
    \end{equation}
    holds for all $\eta\in\mathbb R_{\ge 0}$.    
\end{lem}

\begin{proof}
    To prove the above lemma, we first observe that from the symmetric property of $Q_{sl}^\lambda$, we have
    \[Q_{sl}^\lambda(\eta)\frac{(\lambda+1)}{2\lambda}\le\eta\le Q_{sl}^\lambda(\eta)\frac{(\lambda+1)}{2}.\]
    By taking the upper and lower bound of above inequality, we have that
    \begin{align*}
        \eta\le Q_{sl}^\lambda(\eta)\frac{(\lambda+1)}{2}
        &&\Leftrightarrow&& Q_{sl}^\lambda(\eta) &\ge \left(\frac{2}{\lambda+1}\right)\eta,
    \end{align*}
    and
    \begin{align*}
        \eta\ge Q_{sl}^\lambda(\eta)\frac{(\lambda+1)}{2\lambda}
        &&\Leftrightarrow&& Q_{sl}^\lambda(\eta) &\le \left(\frac{2\lambda}{\lambda+1}\right)\eta,
    \end{align*}
    respectively. By combining both the upper and lower bound of $Q_{sl}^\lambda$, we get
    \[\left(\frac{2}{\lambda+1}\right)\eta\le Q_{sl}^\lambda(\eta)\le \left(\frac{2\lambda}{\lambda+1}\right)\eta.\]
    Finally, 
    the inequality \eqref{ineq:symLogarithmic_quantizer_nnc} holds for all $\eta\geq 0$.
\end{proof}

\begin{prop}\label{cor:nnc_symLog}
    Consider the system $\Sigma_{\rm lin}$ described by \eqref{sys:linear_nnc} satisfying {\rm (A0)} and a discrete set $\mathcal U_l^{\text{ext},\lambda}$ as in \eqref{set:NNLQ} constructed from a finite countable set of unit vectors $\mathcal U$ that together with $\{0\}$
    satisfies {\rm (A2)} so that \eqref{eq:bndVoronoi} holds for some $\delta>0$. Let $\phi_{\mathcal U_l^{\text{ext},\lambda}}$ and $Q_{sl}^\lambda$ be as given in \eqref{eq:nnc_phi_log} and \eqref{eq:logarithmic_symmetric}, respectively, along with a given desired step size $\lambda>1$. Furthermore, assume that {\rm (A1)} holds with $\Psi(y)$ replaced by $-\phi_{\mathcal U_l^{\text{ext},\lambda}}(-y)$ for all $y\in\mathbb R^m$, $k_1=\alpha\left(\frac{2\lambda}{\lambda+1}\right)$ for some $\alpha\in(0,1]$, and $k_2=\left(\frac{2}{\lambda+1}\right)$.
    Then the closed-loop system with $u=\phi_{\mathcal U_l^{\text{ext},\lambda}}(-y)$ is GES.
\end{prop}

The last proposition is similar to the result in Proposition~\ref{cor:nnc_uniform}, with the exception that the assumption (A1) is already satisfied everywhere. Therefore, applying the result in Proposition~\ref{prop:main1} results in the closed-loop system being $0$-GEPS, i.e.\ the closed-loop system is GES.

\section{Illustrative Example}\label{sec:example}

\begin{figure}[hbtp!]
    \centering
    \includegraphics[width=0.75\linewidth]{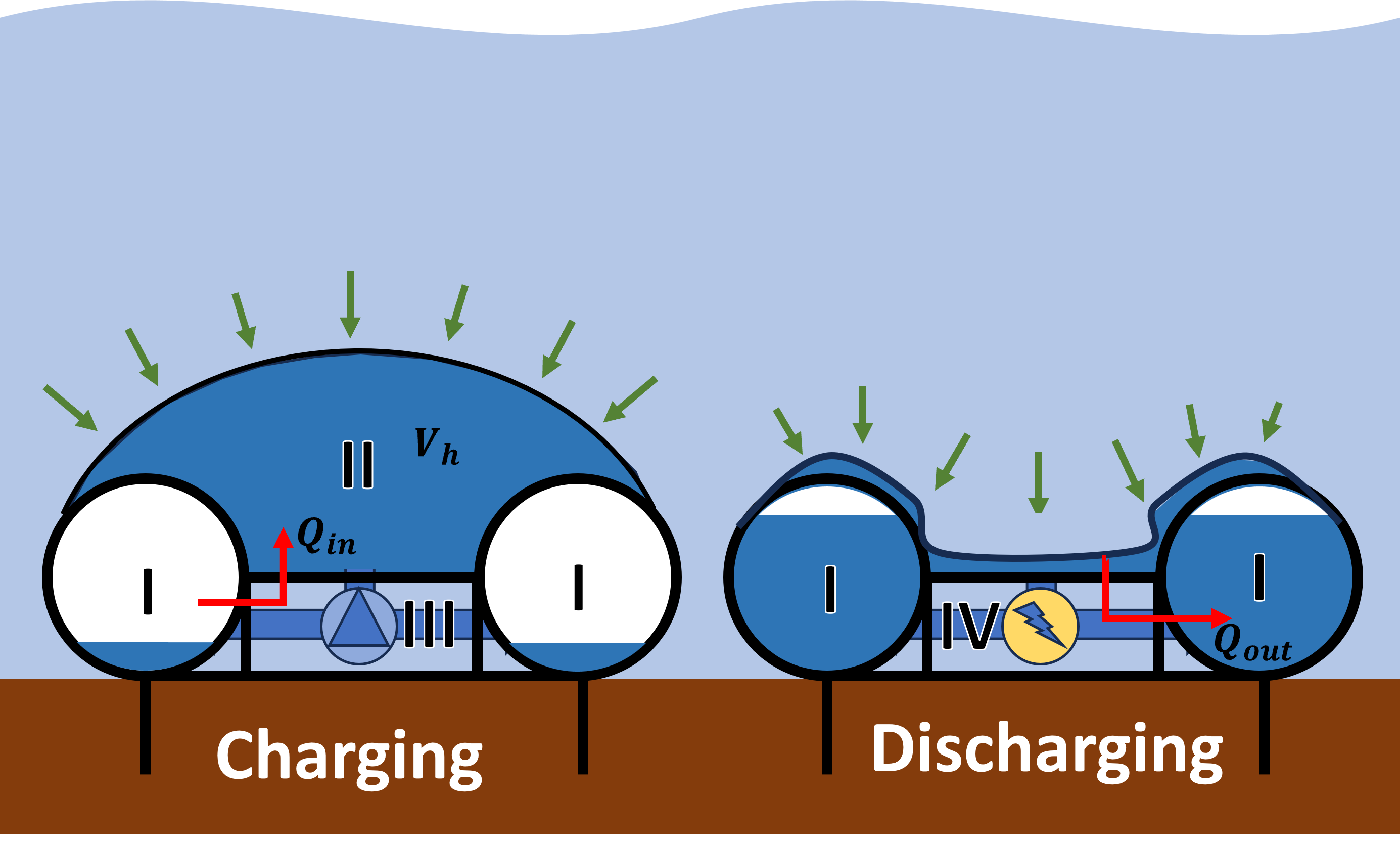}
    \caption{Schematic of a single Ocean Battery system sitting on the seabed with I) low-pressure reservoir (rigid reservoir under atmospheric pressure), II) high-pressure reservoir (flexible reservoir under hydrostatic pressure from underwater environment), III) water pump (blue triangle) to pump water from the low to the high-pressure reservoir, and IV) generator to generate electricity by releasing water from the high- to the low-pressure reservoir. Left picture shows the battery state when charging while the picture on the right side shows the state when discharging.}
    \label{fig:ocean_battery}
\end{figure}

In this section, we present a numerical example to validate our main results for both the uniform and logarithmic nearest-action feedback approaches. We consider an interconnected Ocean Battery system where at least one battery can pump water to the high-pressure reservoir. The Ocean Battery is a novel underwater energy storage system based on the concept of pumped hydro storage \cite{prins2021underwater}. It converts electrical energy supplied from renewable sources such as wind and/or floating solar into stored potential energy which is available naturally due to the presence of hydrostatic pressure from the surrounding environment. For practical applications, the Ocean Battery is installed at a certain depth depending on the requirements of the storage system. For offshore application in sufficiently deep waters, for example, the Ocean Battery is installed on the seabed for the optimal operation; alternatively, the effective depth at which the device operates can be increased by burying the rigid reservoir and machine room into the seabed.

In this example, we consider multiple interconnected Ocean Battery systems whose operations follow the schematic presented in Fig.~\ref{fig:ocean_battery} where $Q_{\text{in}}^i$ and $Q_{\text{out}}^i$ are the total in/outflow of water in the $i$-th battery, and $V_h^i$ is its corresponding water volume in its high-pressure reservoir. We assume that the generator is designed in a way so that the outflow $Q_{\text{out}}^i$ remains constant all the time. A single Ocean Battery system can be described by $\dot V^i = Q^i$ where $V^i=V_h^i-V_h^{i,\rm ref}$ with $V_h^{i,\rm ref}$ be the 
% target
{\startmodif reference}
volume of high-pressure fluid in the $i$-th battery, and $Q^i$ represents the total inflow to the $i$-th battery, i.e.\ $Q^i = Q_{\text{in}}^i-Q_{\text{out}}^i$.

For numerical simulation purposes, we consider the system $\Sigma_{\rm lin}$ in \eqref{sys:linear_nnc} with {\startmodif the state $V=\operatorname{col}\{V^1,\ldots,V^4\}$, a weighted-averaged water volume in the high-pressure reservoirs as the output $y=CV$, water inflow to the high-pressure reservoirs in the ocean batteries which are connected to external sources as the input, and}
\begin{align*}
    A =& \bbm{-0.41 & 0.2 & 0 & 0.2\\
            0.2 & -0.41 & 0.2 & 0\\
            0 & 0.2 & -0.41 & 0.2\\
            0.2 & 0 & 0.2 & -0.41},\\
    & B = \bbm{1 & 0\\ 0&0\\0&1\\0&0},
    \quad C = \bbm{\frac{1}{2} & \frac{1}{4} &0 &\frac{1}{4}\\
    0 & \frac{1}{4} &\frac{1}{2} & \frac{1}{4}}.
\end{align*}
In the context of multiple interconnected Ocean Battery systems, the above system's matrices $(A,B,C)$ can be regarded as an interconnection of distributed  storage devices with 4 storage units where only 2 units are connected to the external (renewable) energy sources for charging and they distribute the stored energy between storage systems as described by the off-diagonal elements of $A$. Each diagonal term in $A$ represents an energy loss and energy transfer to its neighboring storage devices.  
The matrix $B$ defines which storage devices that are connected to the energy generation for charging, %in the system that have pumps 
and the matrix $C$ represents the measured weighted-averaged water volume in the high-pressure reservoir of the system. In this example, all pumping units can only be activated proportionally according to the rule set in a countable set of actions $\mathcal U$. Let the %\textit{directional} 
finite set of unit action vectors $\mathcal U$ be given by
\[\mathcal U:=\left\{ \bbm{\sin{\frac{2\pi}{9}}\\\cos{\frac{2\pi}{9}}}, \bbm{\sin{\frac{8\pi}{9}}\\\cos{\frac{8\pi}{9}}}, \bbm{\sin{\frac{14\pi}{9}}\\\cos{\frac{14\pi}{9}}} \right\}\]
which represents the vertices of an equilateral triangle centered at the origin. The elements in the countable set $\mathcal U$ can be seen as the proportion of charge/discharge in the corresponding batteries. For the set $\mathcal U$ above, the value of $\delta$ satisfying \eqref{eq:bndVoronoi} is $\delta=1$. Moreover, the value of $\alpha$, which is exactly the cosine of the largest possible angle between any points in $\mathbb R^2$ to the nearest point in $\mathcal U$, is $\alpha=\frac{1}{2}$.

For the simulation, we set $\lambda=0.5$ for the uniformly extended quantizer and $\lambda=1.1$ for the logarithmically extended quantizer. 
It can be checked that with the given tuple $(A,B,C)$, the respective transfer function $(I+k_2 G){(I+k_1 G)}^{-1}$ with $k_1$ and $k_2$ as stated in Propositions~\ref{cor:nnc_uniform} and \ref{cor:nnc_symLog} is strictly positive real. To be more precise, the values of $k_1$ and $k_2$ are $k_1=\frac{1}{4}$ and $k_2=\frac{3}{2}$ for the uniformly extended quantizer; and $k_1=\frac{11}{21}$ and $k_2=\frac{20}{21}$ for the logarithmically extended quantizer.

The simulation results in Fig.~\ref{fig:unif_mimo} and \ref{fig:symLog_mimo} confirm that, with additional control actions in each unit vector direction, the practical stabilization and global asymptotic stabilization can be achieved exponentially fast compared to only using $\mathcal U$. %the minimum actions (linear). 
In addition, logarithmically-extended actions render the system asymptotically stable. 

\begin{figure}[t!]
    \centering
    \includegraphics[width=0.95\linewidth]{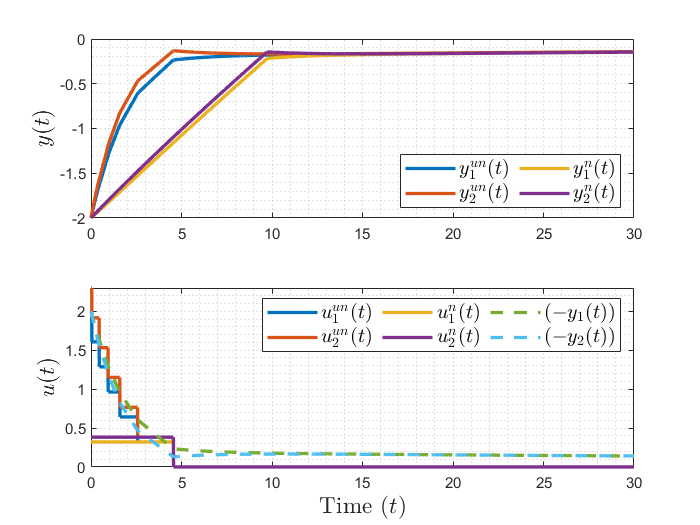}
    \caption{Simulation result with and without additional uniformly distributed quantization levels. The top plot shows the output response where the blue and red lines $y_i^{un}(t),\ i=1,2,$ are the outputs with additional uniformly distributed quantization levels while the yellow and purple lines $y_i^n (t),\ i=1,2,$ are the outputs with only single action in each direction. The bottom plot shows the input signals compared to the negative of the output values.}
    \label{fig:unif_mimo}
\end{figure}

\begin{figure}
    \centering
    \includegraphics[width=0.95\linewidth]{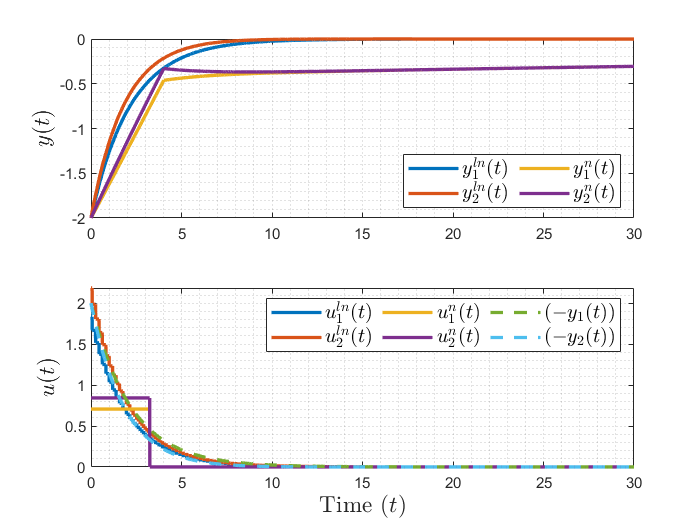}
    \caption{Simulation result with and without additional logarithmically distributed quantization levels. The top plot shows the output response where the blue and red lines $y_i^{un}(t),\ i=1,2,$ are the outputs with additional logarithmically distributed quantization levels while the yellow and purple lines $y_i^n (t),\ i=1,2,$ are the outputs with only single action in each direction. The bottom plot shows the input signals compared to the negative of the output values.}
    \label{fig:symLog_mimo}
\end{figure}

\section{Conclusion}\label{sec:conclusion}
In this work, we propose the use of weak sector condition for MIMO systems and show that the closed-loop system remains {\startmodif ISS}. Moreover, we show that practical stability of the feedback systems with nearest-action input-selection approach can be analyzed using the notion of {\startmodif ISS}. This is achieved by extending the countable action set $\mathcal U$ via uniformly-extended actions or via logarithmically-extended actions. The application of nearest-action control using these extended action sets guarantees that the closed-loop system is exponentially practically stable for the uniform one and globally exponentially stable for the logarithmic one.

\bibliography{bib_all}

\end{document}